\documentclass[12pt]{article} 

\usepackage[utf8]{inputenc} 

\usepackage{geometry} 
\usepackage{graphicx} 

\usepackage{booktabs} 
\usepackage{array} 
\usepackage{paralist} 
\usepackage{verbatim} 
\usepackage{subfig} 

\usepackage{fancyhdr} 
\usepackage{sectsty}
\allsectionsfont{\sffamily\mdseries\upshape} 

\usepackage[nottoc,notlof,notlot]{tocbibind} 
\usepackage[titles,subfigure]{tocloft} 

\usepackage{amsmath}
\usepackage{amsthm}
\newtheorem{example}{Example}
\newtheorem{def1}[example]{Definition}
\newtheorem{def2}[example]{Definition}
\newtheorem{mxp}[example]{Proposition}
\newtheorem{l4}[example]{Lemma}
\newtheorem{c5}[example]{Corollary}

\newtheorem{l6}[example]{Theorem}
\newtheorem{l7}[example]{Theorem}
\newtheorem{l41}[example]{Theorem}

\begin{document}

\title{The peak of the solution of elliptic equations}
\author{Janpou Nee\\
General Education Center, ChienKuo
Technology University\\Changhua, Taiwan\\
jpnee@ctu.edu.tw}
\date{}

\maketitle

%


\begin{abstract}
     A counter example of inheritance of convexity of domain of
positive solution of Dirichlet boundary value problem and the
hot spot problem  that proposed by J. Rauch is given.  The
difficulty of these two problems is that the critical points of
the solutions is not singleton but a level curve.  However, using
Pohozeav identity locally, partial answer to both problems can
be derived.\\
\newline
Keywords: Critical points, inheritance of convexity, Neumann boundary
value problems, hot spot problem.\\
Subjclass[2010] Primary 35J57, 35J61
\end{abstract}

\section{Introduction}

         This article concerns the local behavior of
semi-linear elliptic equation
\begin{equation}\label{or}
     -\Delta u=f(u),\quad x\in D
\end{equation}
with either Neumann boundary condition
\begin{equation}\label{neu}
     \frac{\partial u}{\partial n}|_{\partial D}=0,
\end{equation}
or Dirichlet boundary value problem
\begin{equation}\label{dir}
     u|_{\partial D}=0.
\end{equation}
The solution of equations (\ref{or}) are the steady states of both
heat transferring problem and the standing wave of sound
propagation.  For heat transferring problem, the maximum of the
solution is expected to occur at the boundary.  This problem was
proposed by Rauch \cite{jr}.  The conclusion of this problem is
not necessary true for standing wave problem.  In fact, the counter
example below looks  more like a standing wave solution
than the solution of heat transferring problem.  Perhaps this is the
reason why that there are many counter examples as \cite{bw, b}
and example 1 of this article.  After \cite{bw, b}, this problem then
adjusted to consider over convex domain.  For nonlinear problem,
similar behavior of the solution is called spike layer \cite{lnt}.
Spike layer problem attracts many research attentions, for example
\cite{lin, nt, nit2, ym, jw}.  In this article, the location of the
maximum of the solution of Neumann boundary value and the
uniqueness of local extrema of Dirichlet boundary value problem
are studied.  Example 1 is the counter example of both the hot
spot problem and the inheritance of convexity of positive solution.
It was believed that the positive solution of Dirichlet boundary value
problem of equation (\ref{or}) will inherit the convexity of the domain.
However, counter examples were found by Koreeva \cite{kor},
Cabr\'e and Chanillo \cite{cab} and Hamel et al. \cite{ham} over
convex domain.  Thus, instead of studying the inheritance of convexity
of the solution of elliptic equation, the uniqueness of local maximum
of the solution will be considered in this article.  Although
the condition of the solution that preserve convexity of domain
remains unknown, the condition of positive solution to have a
unique local maximum is obtained.

      At first glance, these two problems seem to be irrelevant
because their boundary conditions are different although they
share the same equation.  However, the behavior of these two
solutions of the contrary proposition of the corresponding problems
s similar.  For Dirichlet boundary condition problem,
the level curves should be a closed curves unless there are more
then one isolated critical points that imply the existence of
saddle point.  On the other hand, the level curves of the solution
of Neumann boundary condition should be orthogonal to the boundary
or the boundary itself is a collection of critical points.  Excluding
the boundary being a collection of critical points, the level curves
of Neumann boundary condition should not form a closed curves
inside the domain unless there exists some other isolated critical
point in the interior of the domain that implies the existence of
saddle point.  Thus, to derive the right condition of these problems,
it is necessary to classify the critical point.  The common
method of classification critical points refers to Hessian matrix.
Unfortunately, the solution of a differential equation is not a definite
function therefore it is not possible to adept the Hessian matrix to it.
To this end, a new testing method is introduced which involves the
following definitions of critical points.

\begin{def1}
An isolated critical point $p$ of function $u$ is a local maximum
(minimum) if for any $\epsilon>0$, there exists a $\delta>0$
such that $\nabla u(x)\cdot(x-p)<0$ ($>$) and
$0<\|\nabla u(x)\|\leq\epsilon$ for
$x\in B_\delta(p)\cap\bar{D}\backslash\{p\}$; otherwise,
$p$ is a saddle point.
\end{def1}

\begin{def2}
     If $s^j_t\subset\mathcal{S}_t$ is a simply connected component
of level curve $\mathcal{S}_t$ and $\nabla u(x)=0$ for all $x\in s^j_t$
then $x$ is called a non-isolated critical point of function $u$.
\end{def2}

     If $p$ is a local maximum of solution $u$ of equation (\ref{or})
then $f(u(p))\geq0$.  To see this, we integrate equation (\ref{or})
over $B_\delta(p)$ and it yields,
     \[-\int_{\partial B_\delta(p)\cap\bar{D}}
     \frac{\partial u}{\partial n}ds
     =\int_{B_\delta(p)\cap\bar{D}}f(u(x))dx.\]
Since $p$ is a local maximum and the outward normal
$n(x)=-\frac{x-p}{\delta}$ along $\partial B_{\delta}(p)$,
it gives $\frac{\partial u}{\partial n}\leq0$.  If
$p\in\partial D$ then $\frac{\partial u}{\partial n}=0$.
Notice that, $\delta$ is arbitrary, therefore
$f(u(p))\geq0$ and hence the following conclusion is derived.
\begin{mxp}
     If $p$ is a local maximum (minimum) of the solution $u$ of
(\ref{or}) then $f(u(p))\geq0$ $(f(u(p))\leq0)$.
\end{mxp}

     Notice that definition 1 and 2 can be applied to classify the
degenerate critical point such as $(0, 0)$ of $x^3-y^3$ or the
critical points at the boundary of the domain.

     Before discussing the main issues, a counter example will
be given first.  The counter example is constructed over disks
because disk posses both convex and symmetry properties.
Thus if the counter example of these problems were found
over disk then it seems that other kinds of necessary conditions
are needed so that these conjectures hold.  In fact, Lemma 5
and the counter example below seem to suggest that the
nonlinearity of reaction term is the key to solve these two
problems.  For instance, the power $m$ of the nonlinear
reaction term of counter example satisfies $m<1$; however,
the condition of Lemma 5 (Pohozeav's identity) says that if
$f(u)=u^m$ then $m>1$.  On the other hand, if $m=1$
then (\ref{or}) turns to eigen-value problem and the location
of hot spot of eigen-function depends on the domain.

\begin{example}
$u(x,y)=\frac{(x^2+y^2)^2}{4}-(x^2+y^2)^{\frac{3}{2}}+(x^2+y^2)$,
which is the symmetry solution over domain $D=B_{r_0}(0)$ of
equation (\ref{sym}) below, where $B_r(p)$ is the ball of radius $r$
centered at point $p$ where
\end{example}
\begin{equation}\label{sym}
     -\Delta u
     =3-\sqrt{1+2\sqrt{u}}+8\sqrt{u},\quad x\in D,
\end{equation}
and its graph is as follows:
\begin{center}
\includegraphics[scale=0.25]{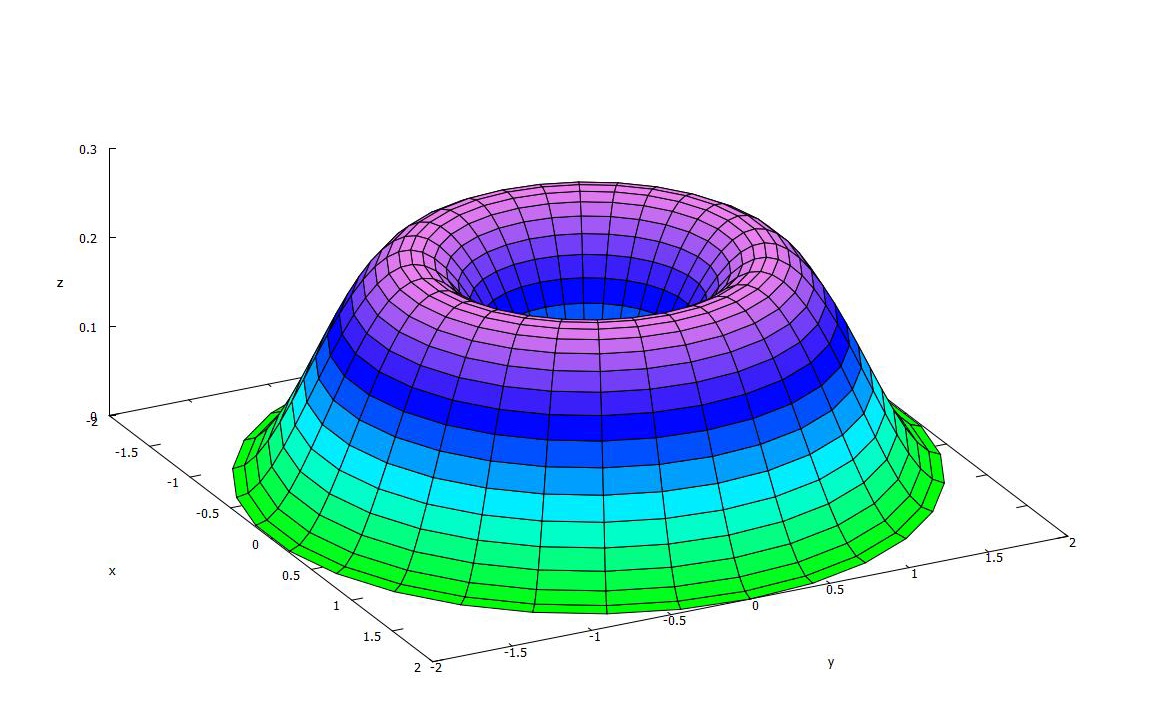}
\end{center}

     If domain $D=B_1(0)$ then $u$ satisfies Neumann boundary
condition and all points of the boundary are maximum of $u$ which
does agree with hot spot conjecture.  However, if $D=B_2(0)$ then
$u$ satisfies both Neumann and Dirichlet boundary condition;
moreover, $u$ is positive but not convex and still $\nabla u=0$
along $B_1(0)$.  Thus $u$ conflicts with both hot spot of heat
transferring on convex domain and the inheritance of convexity
of positive solution.

     In this article, we use $\mathcal{S}_t$ to denote the level
curves of a function.  However, it is possible that the level curves
contains more than one components.  Therefore, we use
$s^j_t\subset\mathcal{S}_t=\cup_{j=1}^ks^j_t$ to denote the
simply connected component of $\mathcal{S}_t$.  Here $s^j_t$
may be a singleton $\{x_0\}$, if $x_0$ is a local extrema.

\section{Main results}

      Throughout this article we assume that $D$ is a smooth
open bounded simply connected convex domain of $R^2$
satisfying interior spherical condition and $f$ is monotone
with respect to $u$ or $\frac{df}{du}>0$.  As usual, we
say that $D$ is smooth if $\partial D$ is smooth.  $\bar{D}$
is the closure of domain $D$ and $\stackrel{\circ}{D}$ is the
interior of $D$.  We shall note that all proofs of this article are
dimensionless; therefore, all results are expected to be true for
higher dimension $N>2$.

     In the proof below, it involves with the critical point at
boundary; therefore, we introduce the following notations
     \[\mathcal{B}_\delta(p)=B_\delta(p)\cap D,\]
     \[\mathcal{B}^+_\delta(p)
     =\{x\in\mathcal{B}_\delta(p):\nabla u(x)\cdot(x-p)>0\},\]
and
     \[\mathcal{B}^-_\delta(p)
     =\{x\in\mathcal{B}_\delta(p):\nabla u(x)\cdot(x-p)<0\}.\]
Without ambiguity, $p\in\partial D$ and $\nabla u(p)=0$ if and
only if $\lim_{x\rightarrow p}\nabla u=0$ where $x\in\bar{D}$.
In particular, if $p\in\partial D$ but $\nabla u(x)\cdot(x-p)$
remains constant sign, for all $x\in\bar{D}$, then we still say
that $p$ is a local extrema.

\subsection{Neumann boundary value problem}

     From Neumann boundary condition, we have a natural
constraint:
\begin{equation}\label{4}
     \int_Df(u)dx=0,
\end{equation}
therefore $f(u)$ must change its sign, say at $u_0$.  As
usual, we let $F(u)=\int f(u)du$.  The assumption
$\frac{df}{du}>0$ implies that $F(t)$ concaves upward with
respect to $t$ and hence $u_0$ is the absolute minimum of
$F$.  To explore the behavior of solution $u$ of Neumann
boundary value problem, we denote $D^+=\{x\in D:f(u(x))>0\}$,
$D^-=\{x\in D:f(u(x))<0\}$ and $m=\min_{x\in\bar{D}}u(x)$
and $M=\max_{x\in\bar{D}}u(x)$.

     Most of the results of this article are based on the following
hypothesis:
\begin{equation}\tag{A}
    \mbox{$u\cdot f(u)-2F(u)>0$},
    \quad\mbox{ $\frac{df}{du}>0$}.
\end{equation}

     First, we consider equation (\ref{or}) with Neumann boundary
value problem.
\begin{l4}
     If $u$ is the smooth solution of (\ref{or}) satisfying
hypothesis (A), $F(t)>0$ and if $p\in\bar{D}$ is an isolated
critical point then $p$ is a local extrema.
\end{l4}
\begin{proof}
The lemma will be proved by deriving a contradiction.
Without loss of generality, we assume that $p\in\bar{D}^+$
such that $\nabla u\cdot(x-p)$ changes its sign.

      To derive the results, we apply Phozeav identity locally.
Multiplying $\nabla u(x)\cdot(x-p)$ to equation (\ref{or})
and integrating over $\mathcal{B}^+_\delta(p)$, it yields
\begin{equation}\label{5}
     \begin{array}{rcl}
     \int_{\mathcal{B}^+_\delta(p)}-\Delta u(\nabla u(x)\cdot(x-p))dx
     &=& -\int_{\partial\mathcal{B}^-_\delta(p)}
              \frac{\partial u}{\partial n}(\nabla u(x)\cdot(x-p))ds\\
      &&  +\int_{\mathcal{B}^-_\delta(p)}
                 \nabla u\cdot\nabla(\nabla u(x)\cdot(x-p))dx.
      \end{array}
\end{equation}
Replacing $-\Delta u$ by $f(u)$, the left hand side of
(\ref{5}) yields
\begin{equation}\label{7}
     \int_{\mathcal{B}^+_\delta(p)}-\Delta u(\nabla u(x)\cdot(x-p))dx
     =\int_{\mathcal{B}^+_\delta(p)}f(u)(\nabla u\cdot(x-p))dx,
\end{equation}
where
\begin{equation}\label{71}
     \int_{\mathcal{B}^+_\delta(p)} f(u)\cdot(\nabla u\cdot(x-p))dx
     =\int_{\mathcal{B}^+_\delta(p)}\nabla F(u)\cdot(x-p)dx,
\end{equation}
and
\begin{equation}\label{72}
     \int_{\mathcal{B}^+_\delta(p)}\nabla F(u)\cdot(x-p)dx
     = \int_{\partial\mathcal{B}^+_\delta(p)}F(u)\cdot((x-p)\cdot n)ds
     -2\int_{\mathcal{B}^+_\delta(p)}F(u)dx.
\end{equation}
Calculating the right hand side of (\ref{5}) and by $D\subset R^2$
it yields
\begin{equation}\label{8}
     \int_{\mathcal{B}^+_\delta(p)}\nabla u\cdot\nabla(\nabla u\cdot(x-p))dx
     =\int_{\partial\mathcal{B}^+_\delta(p)}\frac{|\nabla u|^2}{2}((x-p)\cdot n)ds
     -\int_{\mathcal{B}^+_\delta(p)}\|\nabla u\|^2 dx.
\end{equation}
Replacing $\int_{\mathcal{B}^+_\delta(p)}\|\nabla u\|^2dx$ by
$\int_{\partial\mathcal{B}^+_\delta(p)}\frac{\partial u}{\partial n}uds
+\int_{\mathcal{B}^+_\delta(p)}f(u)udx$ and then adding all together,
we get
\begin{equation}\label{9}
     \begin{array}{rcl}
      0&<&   \int_{\mathcal{B}^+_\delta(p)}-2F(u)+f(u)udx\\
        &=&   \int_{\partial\mathcal{B}^+_\delta(p)}
                  (\frac{\|\nabla u\|^2}{2}-F(u))((x-p)\cdot n)
                 -\frac{\partial u}{\partial n}(\nabla u(x)\cdot(x-p)+u)ds.
     \end{array}
\end{equation}
Let
$\partial\mathcal{B}^+_\delta(p)=\mathcal{N}\cup\mathcal{D}\cup B$
where
$\mathcal{N}=\{x\in\mathcal{B}^+_\delta(p)|\nabla u(x)\cdot(x-p)=0\}$,
$\mathcal{D}=\mathcal{B}^+_\delta(p)\cap\partial D$ and
$B=\partial B_\delta(p)\cap\partial\mathcal{B}^+_\delta(p)$.
If $p\in\partial D$ then $\mathcal{D}\neq\emptyset$ otherwise
it is an empty set.

     Along $\mathcal{N}$, $x-p$ is parallel to the tangent of
the curve therefore $(x-p)\cdot n=0$.  On the other hand,
$\nabla u\cdot(x-p)\geq0$ over $\mathcal{B}^+_\delta(p)$
therefore $\frac{\partial u}{\partial n}\geq0$.  Along $B$,
$x-p\cdot n=\delta$ and $\nabla u\cdot(x-p)\geq0$
therefore $\frac{\partial u}{\partial n}\geq0$.  Along $D$,
$\frac{\partial u}{\partial n}=0$ ($u=0$ for Dirichlet boundary
value problem) and $(x-p)\cdot n=\delta$.  Since $F(u(p))>0$,
(\ref{9}) yields
\begin{equation}\label{10}
     \begin{array}{rcl}
     0&<& \int_{\mathcal{B}^+_\delta(p)}-2F(u)+f(u)udx\\
       & = & \int_{\partial\mathcal{B}^+_\delta(p)}
                 -F(u)((x-p)\cdot n)-\frac{\partial u}{\partial n}u
                 -\frac{\delta\|\nabla u\|^2}{2}ds\\
       & = & \int_{B}-\frac{\delta\|\nabla u\|^2}{2}
                 -\delta F(u)-\frac{\partial u}{\partial n}uds
                -\int_{\mathcal{N}}\frac{\partial u}{\partial n}udx,
     \end{array}    
\end{equation}
Every terms on the right hand side of equation (\ref{10}) are
negative, a contradiction.  Therefore $p$ must be a local extrema.
The proof is completed.
\end{proof}

     From the conclusion of Lemma 5 and implicit function
theorem, the level curves of $u$ are either the union of
disjoint simply connected smooth curves or singletons.  Thus
we have the following conclusion.

\begin{c5}
     If hypothesis (A) holds, $F>0$ and if
$\mathcal{S}_t=\cup_{j=1}^ks^j_t$ then
for all $t$ $s^i_t\cap s^j_t=\emptyset$.
\end{c5}

\begin{l6}
     If hypothesis (A) holds, $F>0$, and if
$\mathcal{S}_{u_0}=s_{u_0}$ contains only one component
then $\mathcal{S}_{u_0}\cap\partial D\neq\emptyset$,
$\nabla u\neq0$ along $\mathcal{S}_{u_0}$ and there exists
a unique $p_\pm\in\partial D^\pm$ such that
$u(p_+)=\max_{x\in\partial D^+}u(x)$ and
$u(p_-)=\min_{x\in\partial D^-}u(x)$, respectively.
\end{l6}
\begin{proof}
     If on the contrary $\mathcal{S}_{u_0}\cap\partial D=\emptyset$
then the sign of $f(u)$ along $\partial D$ remains constant.  If
$u|_{\partial D}=C$ then $\nabla u(x)\cdot T(x)=0$ where $T(x)$
is the unit tangent vector at $x$ along $\partial D$.
Thus $\nabla u(x)=0$ which contradicts Lemma 5.  Hence $u$
cannot be a constant along $\partial D$.  Let
$u(p)=\max_{x\in\partial D}u(x)$ and $u(q)=\min_{x\in\partial D}u(x)$
then by Proposition 4, $f(u(q))\leq0$, a contradiction.  Thus
$\mathcal{S}_{u_0}\cap\partial D=\emptyset$.  By Lemma 5,
$\nabla u\neq0$ along $\mathcal{S}_{u_0}$.

     To prove the uniqueness of local maximum along $\partial  D^+$,
we let $p_i\in\partial D^+$ such that
$u(p_i)=\max_{x\in\partial D^+}u(x)$.  Let $\xi\subset\partial D^+$
be the arc containing all the points lie in between $p_i$.
By mean value theorem, there exists at least a critical point $p_0$
lies in between $p_i$.  By Lemma 5, $p_0$ cannot be a saddle point.
Therefore $p_0$ is a local minimum which contradicts proposition 4.
\end{proof}

     If the condition $F(u)>0$ may relax to $F(u)\geq0$ but
with assumption $F(t)=0$ only when $t=u_0$, then Lemma
5 still holds.  Thus $u$ contains interior maximum provided
that $\mathcal{S}_{u_0}$ contains more than one component.

\subsection{Dirichlet boundary value problem}

     Lemma 5 is a local property of the solution of equation (\ref{or})
therefore it remains true for Dirichlet boundary value problem.  With
the conclusion of Lemma 5 and mean value theorem, we may derive
that the positive solution of (\ref{or}) has a unique local maximum
provided that the domain is convex.  From the proof of Lemma 5, we
see that the necessary condition of it is that $F>0$.  The positiveness
of the solution and the conditions $\frac{df}{du}>0$ and $f(t)>0$
for $t>0$ imply that $F(u)>0$ over $\stackrel{\circ}{D}$.  Thus the
following conclusion holds.

\begin{l41}
     If $u$ is the smooth positive solution of (\ref{or}) satisfying
Dirichlet boundary condition with hypothesis (A), $f(0)\geq0$
then $u$ has a unique local maximum.
\end{l41}
\begin{proof}
     If $p$, $q$ are both local maximum of $u$ then by mean value
theorem there must another critical point $p_0$ which is either a
saddle or a local minimum.  If $p_0\in\stackrel{\circ}{D}$ then
$p_0$ cannot be a local minimum because $f(u)>0$ which contradicts
proposition 4.  $p_0$ cannot be a saddle because that will
contradict Lemma 5.  Next, if $p_0\in\partial D$ and if it is a
saddle or local minimum then there is a subset
$B^+\subset B_\delta(p)\cap D$ such that $\nabla u\cdot (x-p_0)>0$,
if $x\neq p_0$, which contradicts Lemma 5.  Thus the interior local
extrema is unique.  The proof is completed. 
\end{proof}

{\bf Remark 9.} The assumptions $f(u)u>2F(u)$ and $F(u)>0$
of hypothesis (A) indicate that if $f(u)=u^p$ then $p>1$ which
coincides with counter example 1.  The condition $\frac{df}{du}>0$
does fit the first non-constant eigenfunction of Laplacian with Neumann
boundary condition.  However, excluding the constant eigenfunction,
the second non-constant eigenfunction $cos(x)cos(y)$ on
$[0, 2\pi] \times[0, 2\pi]$ has an interior critical point and level
curve $\mathcal{S}_0$ contains two components.  Therefore, the
location of the local extrema seems not only depends on the
convexity of the domain but the structure of the level curve of $u_0$
as well.

\bibliographystyle{amsplain}

\end{document}